# An Upwind Generalized Finite Difference Method for Meshless Solution of Two-phase Porous Flow Equations

Xiang Rao[1, 2, 3]*, Yina Liu[2, 3], Hui Zhao[2, 3]*

[1]Cooperative Innovation Center of Unconventional Oil and Gas (Ministry of Education & Hubei Province), Yangtze University, Wuhan, 430100, China

[2]Key Laboratory of Drilling and Production Engineering for Oil and Gas, Hubei Province, Wuhan 430100, Hubei, China

[3]School of Petroleum Engineering, Yangtze University, Wuhan 430100, China

*Corresponding author: Xiang Rao (raoxiang0103@163.com, raoxiang@yangtzeu.edu.cn), Hui Zhao (zhaohui@yangtzeu.edu.cn)

Abstract: This paper makes the first attempt to apply newly developed upwind GFDM for the meshless solution of two-phase porous flow equations. In the presented method, node cloud is used to flexibly discretize the computational domain, instead of complicated mesh generation. Combining with moving least square approximation and local Taylor expansion, spatial derivatives of oil-phase pressure at a node are approximated by generalized difference operators in the local influence domain of the node. By introducing the first-order upwind scheme of phase relative permeability, and combining the discrete boundary conditions, fully-implicit GFDM-based nonlinear discrete equations of the immiscible two-phase porous flow are obtained and solved by the nonlinear solver based on the Newton iteration method with the automatic differentiation, to avoid the additional computational cost and possible computational instability caused by sequentially coupled scheme. Two numerical examples are implemented to test the computational performances of the presented method. Detailed error analysis finds the two sources of the calculation error, roughly studies the convergence order thus find that the low-order error of GFDM makes the convergence order of GFDM lower than that of FDM when node spacing is small, and points out the significant effect of the symmetry or uniformity of the node collocation in the node influence domain on the accuracy of generalized difference operators, and the radius of the node influence domain should be small to achieve high calculation accuracy, which is a significant difference between the studied hyperbolic two-phase porous flow problem and the elliptic problems when GFDM is applied. In all, the upwind GFDM with the fully implicit nonlinear solver and related analysis about computational performances given in this work may provide a critical reference for developing a general-purpose meshless numerical simulator for porous flow problems.

1. Introduction

The meshless generalized finite difference method (GFDM) is a meshless method [1,2] developed in recent years. Based on the Taylor series expansion of multivariate functions in the subdomain of each node (i.e. the node influence domain) and the weighted least squares approximation, the first-order and the second-order spatial derivatives of the unknown function in the differential equation are expressed as difference schemes of nodal values in the subdomain, which overcomes the dependence of the traditional finite difference method on the mesh division. At present, this method has developed rapidly and is widely used to solve various scientific and engineering problems, including shallow water equation [3], high-order partial differential equation [4], transient heat conduction analysis [5], stress analysis [6], water wave interaction [7], inverse heat source problem [8], seismic wave propagation problem [9], coupled thermoelasticity problem [10-12], stochastic analysis of groundwater flow [13], And some typical differential equation problems (e.g., unsteady Burger's equation [14-15], nonlinear convection-diffusion equation [16], fractional diffusion equation [17]). GFDM only needs to allocate nodes in the calculation domain to realize the effective solution of the differential equations, which can avoid the complicated meshing of the calculation domain with complex geometry in mesh-based methods and a large number of numerical integrals required in some other meshless methods. Benito et al. [18] developed an h-adaptive method in GFDM.

For flow problems, it is often necessary to take an upwind scheme when handling some physical parameters in governing equations. In the meshless method, the upstream scheme is generally realized by modifying the influence domain of the node, including the upwind influence domain method [19] that moves the central node to the upstream direction and the partial influence domain method [20] that more upstream nodes are included in the influence domain of the central node. However, because the actual flow field may be very complex, it is difficult to form a stable upwind effect by modifying the node influence domain, and satisfactory calculation accuracy is difficult to be well guaranteed. Sridar and Balakrishnan [21] presented an upwind least-squares-based finite difference method for computational fluid dynamics. Saucedo-Zendejo et al. [22] adopted GFDM to solve 3D free surface flows in casting mould filling processes with a semi-implicit scheme. Michel et al. [23] applied GFDM to simulate solution mining processes on microscopic and macroscopic scales. Shao et al. [24] used GFDM to solve Stokes interface problems. However, the subsurface porous flow problems have not been studied by GFDM.

Rao et al. [25] applied an upwind GFDM to model single-phase heat and mass transfer in porous media with a sequential coupled scheme and proposed a GFDM based method to solve the convection-diffusion equation with high accuracy, which also implied the great potential of the upwind GFDM in porous flow problems. Immiscible two-phase flow is a basic scenario of porous flow, such as oil-water two-phase flow in an underground reservoir [26, 27, 28]. The relative phase permeability in the two-phase porous flow equations is a typical physical parameter that needs to take the upwind scheme. To verify the application of the upwind GFDM in porous multiphase flow problems, this paper for the first time studies the computational performances of the upwind GFDM in two-phase porous flow, and solve the equations by using the fully implicit scheme and the Newton method based nonlinear solver, so as to eliminate the additional computational cost and possible computational instability in the sequentially coupled scheme. Numerical examples in this paper validate the good computational performances of the upwind GFDM for two-phase porous flow problems. In addition, this paper conducts the detailed error analysis, and to our knowledge for the first time points out: (i) due to lower-order error of GFDM resulted by the weighted least squares treatment compared with the truncation error of FDM, when the node spacing is large, the convergence order is nearly the same as that of FDM, but when the node spacing is small, low-order error of GFDM will be strong such that the convergence order of GFDM is lower than that of FDM. (ii) the two-phase porous flow problem the radius of the node influence domain in GFDM is required small to ensure computational accuracy, but for elliptic problems, there is no significant requirement for the radius of the node influence area in GFDM.

This paper is structured as follows: Section 2.1 introduces the basic governing equations of the two-phase porous flow problem. Section 2.2 presents the GFDM theory. Section 2.3 derives the upwind GFDM-based fully-implicit discrete schemes of the governing equations and boundary conditions. Two numerical examples are implemented in Section 3. Section 4 conducts a detailed error analysis, in which, Section 4.1 the effect of the allocation uniformity in the node influence domain on the calculation accuracy of the two-phase porous flow hyperbolic problem. Section 4.2 points out the large radius of the node influence domain will make the discontinuous water drive front unclear to induce calculation errors. Section 4.3 analysis the differences of the hyperbolic problem and the elliptic problems on the requirement of the radius of the node influence domain in GFDM. Section 4.4 roughly compares the convergence order of the upwind GFDM with that of the upwind FDM. The conclusions come in Section 5.

2. Methodology
2.1 Coupled governing equations in two-phase porous flow

Taking the immiscible oil-water two-phase flow in porous media as an example, the governing equations of the two-phase flow include the mass conservation equations of water phase and oil phase, and the auxiliary equations representing the relationship between the physical quantities in the mass conservation equation, respectively:

*Mass conservation equations of oil phase*:

$$\alpha \nabla \cdot \left( \frac{k k_{ro}}{\mu_o} \nabla p_o \right) + q_o = \frac{\partial (\phi S_o)}{\partial t} \quad (1)$$

where $k$ is the permeability, mD (i.e. $10^{-15}$ m$^2$); $k_{ro} = k_{ro}(S_w)$ is the relative permeability of the oil phase, which is a function of water saturation $S_w$; $\mu_o$ is the oil-phase viscosity, mPa·s; $p_o$ is the oil-phase pressure, MPa; $q_o$ is the source or sink item of the oil phase, 1/day; $t$ is time, day; $\phi = \phi(p)$ is the reservoir porosity which is a function of pressure, fraction; $S_o$ is the oil saturation; The unit of length is meter. $\alpha$ is the constant when each physical quantity takes the above unit, equal to 0.0864.

*Mass conservation equations of water phase*:

$$\alpha \nabla \cdot \left( \frac{k k_{rw}}{\mu_w} \nabla p_w \right) + q_w = \frac{\partial (\phi S_w)}{\partial t} \quad (2)$$

where, $k_{rw} = k_{rw}(S_w)$ is the relative permeability of the water phase, which is a function of the water saturation;

$\mu_w$ is the water-phase viscosity, mPa·s; $p_w$ is the water-phase pressure, MPa; $q_w$ is the source or sink item of the water phase, 1/day;

*Auxiliary equations*:

The basic variables in the above two equations are generally selected as oil phase pressure $p_o$ and water saturation $S_w$. The functional relationship between other physical quantities and the two basic variables is as follows:

$$S_o = 1 - S_w, \quad p_w = p_o - P_c(S_w), \quad \phi = \phi(p_o) = \phi_0 + C_r(p_o - p_o^0), \quad k_{rw} = k_{rw}(S_w), \quad k_{ro} = k_{ro}(S_w) \tag{3}$$

where, $C_t = C_t(S_w)$ is the comprehensive compressibility coefficient, 1/MPa; $p_o^0$ is the original oil phase pressure and $\phi_0$ is the porosity at $p_o^0$; $\phi = \phi(p)$ is the porosity which is a function of the oil-phase pressure; $k_{rw} = k_{rw}(S_w)$ and $k_{ro} = k_{ro}(S_w)$ are the relative permeabilities of the water phase and the oil phase which are functions of the water saturation. Therefore, Eq. (1) and Eq. (2) constitute a nonlinear system about $p_o$ and $S_w$. In addition, the capillary force $P_c(S_w)$ is often ignored or assumed small in the study of two-phase flow, so the two-phase flow problem studied in this paper can be seen a hyperbolic problem [29]. For the sake of brevity, suppose the capillary force $P_c(S_w) = 0$ in the following text.

2.2. Meshless upwind GFDM for the two-phase porous flow problems

GFDM is a meshless collocation method developed rapidly in recent years. By using local Taylor expansion and moving least squares in the node influence domain, the each-order spatial derivative of unknown variables at the node can be approximated in a difference scheme, which is a linear combination of the variable values of each node in the influence domain of the node.

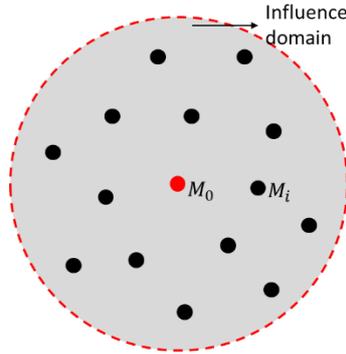

Fig. 1 A simple sketch of a node influence domain and its containing nodes

As shown in Fig. 1, for a node $M_0 = (x_0, y_0)$ in the calculation domain, suppose another $n$ nodes are within the influence domain of the node $M_0$, and denoted as $\{M_1, M_2, M_3, \cdots, M_n\}$, where $M_i = (x_i, y_i)$. Take Taylor expansion of $\{u(M_i), i = 1, \cdots n\}$ at the node $M_0 = (x_0, y_0)$ to obtain:

$$u|_{M_i} = u|_{M_0} + \Delta x_i \, u_x|_{M_0} + \Delta y_i \, u_y|_{M_0} + \frac{1}{2}\left((\Delta x_i)^2 u_{xx}|_{M_0} + 2\Delta x_i \Delta y_i \, u_{xy}|_{M_0} + (\Delta y_i)^2 u_{yy}|_{M_0}\right) + o(\Delta x_i^2, \Delta y_i^2) \tag{4}$$

where $\Delta x_i = x_0 - x_i$, $\Delta y_i = y_0 - y_i$.

It can be seen that the closer $M_i = (x_i, y_i)$ is to the central node $M_0 = (x_0, y_0)$, the right side of Eq. (4) after removing the remainder should be more valid. Therefore, define the weighted error function in Eq. (5), in which the closer the distance, the greater the weight.

$$B(\mathbf{D}_u) = \sum_{j=1}^n \left[\left(u_0 - u_j + \Delta x_j \, u_x|_{M_0} + \Delta y_j \, u_y|_{M_0} + \frac{1}{2}(\Delta x_j)^2 u_{xx}|_{M_0} + \frac{1}{2}(\Delta y_j)^2 u_{yy}|_{M_0} + \Delta x_j \Delta y_j \, u_{xy}|_{M_0}\right)\omega_j\right]^2 \tag{5}$$

where, $\mathbf{D}_u = \left(u_x|_{M_0}, u_y|_{M_0}, u_{xx}|_{M_0}, u_{yy}|_{M_0}, u_{xy}|_{M_0}\right)^T$, $\omega_j = \omega(\Delta x_j, \Delta y_j)$ is the weight function, and the quartic spline function in Eq. (6) is generally selected to calculate the weight in Eq. (5):

$$\omega_j = \begin{cases} 1 - 6\left(\dfrac{r_j}{r_e}\right)^2 + 8\left(\dfrac{r_j}{r_e}\right)^3 - 3\left(\dfrac{r_j}{r_e}\right)^4 & r_j \leq r_e \\ 0 & r_j > r_e \end{cases} \tag{6}$$

where $r_j = \sqrt{\Delta x_i^2 + \Delta y_i^2}$ is the Euclidean distance from node $M_j$ to node $M_0$, $r_e$ is the radius of the node influence domain. In Section 4, we will point out that the radius of the node influence domain should be as small as possible to achieve high calculation accuracy. This is a significant difference between the studied hyperbolic two-phase porous flow problem and the elliptic equation previously studied by GFDM, because the radius of the node influence domain in a large number of studies using GFDM to solve the elliptic equation is often large.

Based on the moving least squares method, the minimum value of $B(\mathbf{D}_u)$ is required, then the partial derivatives of $B(\mathbf{D}_u)$ with respect to each component in the independent variable $\mathbf{D}_u$ should be equal to zero, these are,

$$\frac{\partial B(\mathbf{D}_u)}{\partial u_x}\bigg|_{M_0}=0, \frac{\partial B(\mathbf{D}_u)}{\partial u_y}\bigg|_{M_0}=0, \frac{\partial B(\mathbf{D}_u)}{\partial u_{xx}}\bigg|_{M_0}=0, \frac{\partial B(\mathbf{D}_u)}{\partial u_{yy}}\bigg|_{M_0}=0, \frac{\partial B(\mathbf{D}_u)}{\partial u_{xy}}\bigg|_{M_0}=0 \tag{7}$$

Eq. (7) is sorted into linear equations as follows:

$$\mathbf{A}\mathbf{D}_u = \mathbf{b} \tag{8}$$

where $\mathbf{A}=\mathbf{L}^T\boldsymbol{\omega}\mathbf{L}$, $\mathbf{b}=\mathbf{L}^T\boldsymbol{\omega}\mathbf{U}$, $\mathbf{L}=\left(\mathbf{L}_1^T,\mathbf{L}_2^T,\cdots,\mathbf{L}_n^T\right)^T$, $\mathbf{L}_i=\left(\Delta x_i,\Delta y_i,\frac{\Delta x_i^2}{2},\frac{\Delta y_i^2}{2},\Delta x_i\Delta y_i\right)$, $\boldsymbol{\omega}=diag\left(\omega_1^2,\omega_2^2,\cdots,\omega_n^2\right)$, $\mathbf{U}=\left(u_1-u_0,u_2-u_0,\cdots,u_n-u_0\right)^T$.

Then it is obtained that,

$$\mathbf{D}_u=\left(u_{x0},u_{y0},u_{xx0},u_{yy0},u_{xy0}\right)^T=\mathbf{A}^{-1}\mathbf{b}=\mathbf{A}^{-1}\mathbf{L}^T\boldsymbol{\omega}\mathbf{U}=\mathbf{E}\mathbf{U} \tag{9}$$

where $\mathbf{E}=\mathbf{A}^{-1}\mathbf{L}^T\boldsymbol{\omega}$。

Thus, the approximated difference scheme of the first-order and second-order spatial derivatives at the central node are:

$$u_x\big|_{M_0}=\sum_{j=1}^{n}e_{1j}\left(u\big|_{M_j}-u\big|_{M_0}\right), \quad u_y\big|_{M_0}=\sum_{j=1}^{n+1}e_{2j}\left(u\big|_{M_j}-u\big|_{M_0}\right), \quad u_{xx}\big|_{M_0}=\sum_{j=1}^{n+1}e_{3j}\left(u\big|_{M_j}-u\big|_{M_0}\right),$$
$$u_{yy}\big|_{M_0}=\sum_{j=1}^{n+1}e_{4j}\left(u\big|_{M_j}-u\big|_{M_0}\right), \quad u_{xy}\big|_{M_0}=\sum_{j=1}^{n+1}e_{5j}\left(u\big|_{M_j}-u\big|_{M_0}\right) \tag{10}$$

where $e_{ij}$ is the element of matrix $\mathbf{E}$.

Therefore, GFDM can obtain the approximation of the first-order and the second-order spatial derivatives at the central node when only need to know the coordinates of nodes in the influence domain of the central node, and can iteratively calculate the difference scheme of higher-order spatial derivatives at the central node, so that the meshless discretization of the governing differential equations can be realized without complicated mesh generation. In the discretization of the computational domain with complex geometry, meshless GFDM has a significant advantage over traditional mesh-based methods.

In practice, commonly-used mesh methods, such as the finite element method (FEM), actually generate nodes (i.e. mesh vertices) when generating meshes. Therefore, the mesh generation technology in these mesh-based methods can be directly used to generate the point (node) cloud in the calculation domain. The numerical example in Section 3.1 uses the Cartesian collocation in the traditional FDM constructed for the rectangular calculation domain. The example in Section 3.2 adopts the nodes corresponding to Delaunay triangulation in FEM constructed for the calculation domain with irregular polygonal boundary geometry. In addition to using mesh generation technology to generate corresponding point clouds in the calculation domain, Milewski [30] pointed out that the Liszka-type nodes generator can be used to generate irregular point clouds, which was proposed by Liszka [31] can take full advantage of the irregular calculation domain. Michel et al. [23] applied a meshfree advancing front technique [32] to generate the initial point cloud. Therefore, the technology of generating point cloud in the calculation domain is more diverse than that of generating mesh in the calculation domain. Moreover, Milewski [30] also pointed out that the generated point cloud is not limited by any structure, but the mesh needs to determine a richer topological structure based on the point cloud, including which two nodes connect to form an edge, which nodes form a mesh, the order of mesh vertices, etc. Therefore, although the generation of the point cloud can borrow the technology of mesh generation, it only borrows the part of generating point cloud in mesh generation technology, without the subsequent tedious determination of richer topology on the point cloud. Therefore, the complexity of generating point cloud is significantly lower than that of generating the mesh in the calculation domain. The above analysis supports a significant advantage of meshless methods which only rely on point cloud compared with mesh-based methods in the discretization of the calculation domain. As shown in Fig. 2 (a), taking a simple circular calculation domain as an example, it can be intuitively explained that when the calculation domain is discretized with a point cloud, as shown in

Fig. 2 (b), we can choose to allocate nodes on concentric arcs with equal radius interval and equal angle interval of 60 °. When the circular domain needs to be meshed, it is obviously much more complex than the point allocation strategy in Fig. 2 (b). Even if the mesh is further formed on the basis of the node collocation in Fig. 2 (b), the criteria for determining which points constitute a mesh need to be given, which is much complicated.

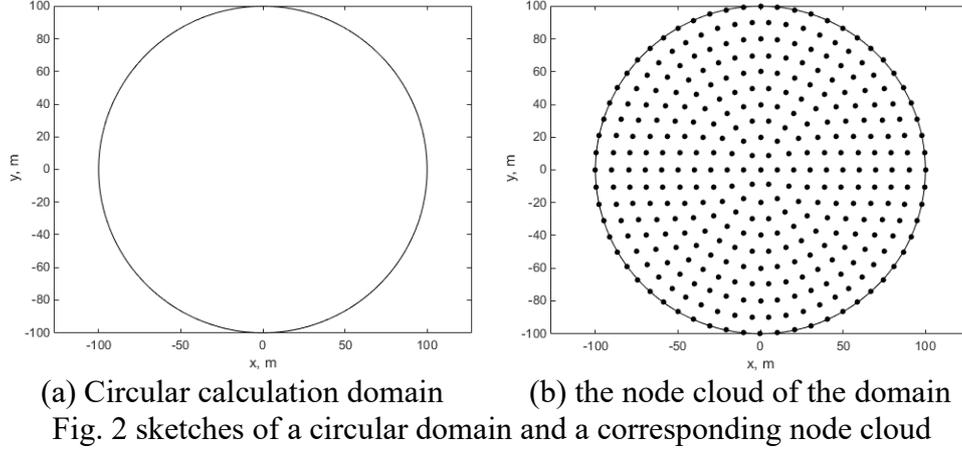

(a) Circular calculation domain    (b) the node cloud of the domain
Fig. 2 sketches of a circular domain and a corresponding node cloud

2.3 Fully implicit GFDM-based discrete schemes

In this section, the fully implicit coupling scheme is used to discretize the governing equations, in which the node pressure and saturation at the time $t+\Delta t$ are solved simultaneously. It may be assumed that the source or sink term is zero, the discrete scheme of oil-phase and water-phase equations based on GFDM is as follows:

$$\alpha \sum_{j=1}^{n_i} \left[ \frac{k_{ij}^{t+\Delta t} k_{ro,ij}^{t+\Delta t}}{\mu_{o,ij}^{t+\Delta t}} \left(e_{3,j}^i + e_{4,j}^i\right) \left(p_{o,(i,j)}^{t+\Delta t} - p_{o,i}^{t+\Delta t}\right) \right] = \frac{1}{\Delta t}\left[\phi_i^{t+\Delta t}\left(1-S_{w,i}^{t+\Delta t}\right) - \phi_i^t\left(1-S_{w,i}^t\right)\right] \tag{11}$$

$$\alpha \sum_{j=1}^{n_i} \left[ \frac{k_{ij}^{t+\Delta t} k_{rw,ij}^{t+\Delta t}}{\mu_{w,ij}^{t+\Delta t}} \left(e_{3,j}^i + e_{4,j}^i\right) \left(p_{o,(i,j)}^{t+\Delta t} - p_{o,i}^{t+\Delta t}\right) \right] = \frac{\left[\phi_i^{t+\Delta t} S_{w,i}^{t+\Delta t} - \phi_i^t S_{w,i}^t\right]}{\Delta t} \tag{12}$$

where $\phi_i^{t+\Delta t} = \phi\left(p_{o,i}^{t+\Delta t}\right)$, $\phi_i^t = \phi\left(p_{o,i}^t\right)$, $P_{c,(i,j)}^{t+\Delta t} = P_c\left(S_{w,(i,j)}^{t+\Delta t}\right)$, $P_{c,i}^{t+\Delta t} = P_c\left(S_{w,i}^{t+\Delta t}\right)$. The permeability $k_{ij}^{t+\Delta t}$ and the viscosity $\mu_{w,ij}^{t+\Delta t}$ adopt harmony average scheme and arithmetic average scheme [26, 33, 34], respectively, these are:

$$k_{ij}^{t+\Delta t} = \frac{2}{1/k_i^{t+\Delta t} + 1/k_j^{t+\Delta t}}, \quad \mu_{o,ij}^{t+\Delta t} = \frac{\mu_{o,i}^{t+\Delta t} + \mu_{o,j}^{t+\Delta t}}{2}, \quad \mu_{w,ij}^{t+\Delta t} = \frac{\mu_{w,i}^{t+\Delta t} + \mu_{w,j}^{t+\Delta t}}{2} \tag{13}$$

Eq. (14) gives the upwind scheme of phase relative permeabilities. It can be seen that the upwind GFDM adopts the first-order upwind scheme without modifying the node influence domain to realize the upwind effect.

$$k_{ro,ij}^{t+\Delta t} = \begin{cases} k_{ro,i}\left(S_{w,i}^{t+\Delta t}\right) & \text{if } p_{o,(i,j)}^{t+\Delta t} \geq p_{o,i}^{t+\Delta t} \\ k_{ro,i}\left(S_{w,i}^{t+\Delta t}\right) & \text{if } p_{o,(i,j)}^{t+\Delta t} < p_{o,i}^{t+\Delta t} \end{cases}, \quad k_{rw,ij}^{t+\Delta t} = \begin{cases} k_{rw}\left(S_{w,(i,j)}^{t+\Delta t}\right) & \text{if } p_{o,(i,j)}^{t+\Delta t} \geq p_{o,i}^{t+\Delta t} \\ k_{rw}\left(S_{w,i}^{t+\Delta t}\right) & \text{if } p_{o,(i,j)}^{t+\Delta t} < p_{o,i}^{t+\Delta t} \end{cases} \tag{14}$$

It should be noted that in this paper, as shown in Fig. 1, the derivative boundary conditions in the calculation domain are treated by adding virtual nodes outside the boundary.

Milewski [30] developed a high-order least squares meshless difference method to improve the calculation accuracy, and pointed out that in the meshless difference method, the local difference approximation accuracy will be reduced due to the low quality of the node collocation in the influence domain of boundary nodes or nodes near the boundary, and briefly introduced that the method of adding virtual nodes outside the boundary can be used to improve the accuracy of the local approximation. In Section 4, we will take Cartesian collocation as an example to analyze in detail the reasons for adding virtual nodes with specific data. In this section, we first point out that as shown in Fig. 3, the derivative boundary conditions are handled by adding a virtual node correspondingly to each boundary node along the external normal vector at the boundary node with a certain distance.

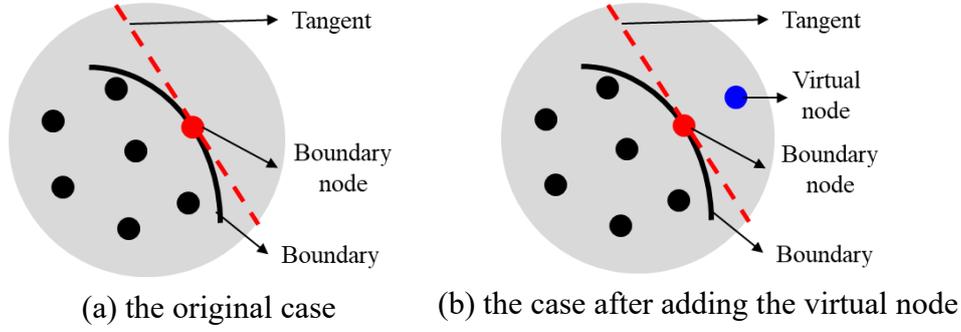

(a) the original case      (b) the case after adding the virtual node

Fig. 3 Sketch of adding virtual nodes [25]

In the calculation domain, suppose there are $n_1$ internal nodes, $n_2$ nodes with Dirichlet boundary conditions, and $n_3$ nodes with derivative boundary conditions.

For a node $A$ with the derivative boundary condition, suppose the sequence number of node $A$ in all nodes be $a$, and the added virtual node corresponding to node $A$ be recorded as node $B$. Since each derivative boundary condition node needs to add a corresponding virtual node, the number of nodes in the whole calculation domain has $n_1+n_2+n_3+n_3$ nodes. Suppose the sequence number of virtual node $B$ in all nodes be $b$, then the equation at node $A$ is no longer the equation corresponding to the boundary condition, but the discrete schemes of the governing equations shown in Eq. (11) and Eq. (12), which are the same as those at each internal node. Suppose the derivative boundary condition of pressure or water saturation is in the form of Eq. (15), the discrete schemes subject to the derivative boundary conditions in Eq. (16) are used as the equations at the virtual node $B$.

$$\left(\alpha u + \beta \frac{\partial u}{\partial \mathbf{n}}\right)\bigg|_A = \gamma, \quad u = p, S_w \tag{15}$$

$$\alpha u_a + \beta \sum_{j}^{n_A} \left(n_x e_{1j} u_{(a,j)} + n_y e_{2j} u_{(a,j)}\right) = \gamma, \quad u = p, S_w \tag{16}$$

where $\alpha$, $\beta$ and $\gamma$ are coefficients. $\mathbf{n}$ is the unit external normal vector at node $A$, and $\mathbf{n} = (n_x, n_y)$. Because node $B$ is within the node influence domain of node $A$, there must be a $j$ that meets $(a, j) = b$.

For a node $C$ with the Dirichlet boundary condition, no corresponding virtual node is added. Therefore, suppose the Dirichlet boundary condition is in the form of Eq. (17) and the sequence number of node $C$ in all nodes be $c$, the equations at node $C$ are just the discrete schemes in Eq. (18) of the Dirichlet boundary condition.

$$u|_C = \eta, \quad u = p, S_w \tag{17}$$

$$u_c = \eta, \quad u = p, S_w \tag{18}$$

where $\eta$ is a coefficient.

Finally, a system of equations composed of $2(n_1+n_2+n_3+n_3)$ equations, including $n_1+n_3$ Eq. (11), $n_1+n_3$ Eq. (12), $2n_2$ Eq. (18)s, $2n_3$ Eq. (16). In this paper, the Newton iteration- and automatic differentiation- based nonlinear solver is used to obtain a stable solution of the fully-implicit nonlinear equations [26, 32, 33]. Then the pressure and water saturation values of all nodes (including $n_1$ inner nodes, $n_2+n_3$ boundary nodes, and $n_3$ virtual nodes) at the time $t+\Delta t$ are solved to avoid the additional calculation cost of small time steps which is adopted to ensure the numerical stability and accuracy in the sequentially coupled scheme.

3. Numerical cases

In this section, a numerical case with a Cartesian node cloud in a rectangular domain and the other case with irregular node clouds in rectangular and polygonal domains are given to test the computational performances of the upwind GFDM.

3.1 A case with a Cartesian node cloud in a rectangular domain

This example selects a regular rectangular computational domain ([0m, 200m] × [0m, 80m]), suppose the left boundary and right boundary (denoted as $\Gamma_1$ and $\Gamma_2$, respectively) meet the first type of boundary conditions, and the upper boundary and lower boundary (denoted as $\Gamma_3$ and $\Gamma_4$, respectively) meet closed boundary conditions, and the source or sink term is zero. Table 1 shows the values of physical parameters, and Eq. (19) summarizes the governing equations with initial and boundary conditions.

$$\alpha \nabla \cdot \left( \frac{kk_{ro}}{\mu_o} \nabla p_o \right) = \frac{\partial (\phi S_o)}{\partial t}, \quad \alpha \nabla \cdot \left( \frac{kk_{rw}}{\mu_w} \nabla p_w \right) = \frac{\partial (\phi S_w)}{\partial t},$$

$$k_{rw} = \left( \frac{S_w - S_{wc}}{1 - S_{or} - S_{wc}} \right)^2, \quad k_{ro} = \left( \frac{1 - S_w - S_{or}}{1 - S_{wc} - S_{or}} \right)^2, \quad p_0 = 10, \quad S_w = 0.8, \qquad (19)$$

$$p|_{\Gamma_1} = 15, \quad p|_{\Gamma_2} = 10, \quad \frac{\partial p}{\partial y}\bigg|_{\Gamma_3} = 0, \quad \frac{\partial p}{\partial y}\bigg|_{\Gamma_4} = 0, \quad S_w|_{\Gamma_1} = 0.8, \quad S_w|_{\Gamma_2} = 0.2, \quad \frac{\partial S_w}{\partial y}\bigg|_{\Gamma_3} = 0, \quad \frac{\partial S_w}{\partial y}\bigg|_{\Gamma_4} = 0$$

In this example, the FDM results obtained by small space steps and small time steps ($\Delta x = 0.1m$, $\Delta y = 0.1m$, $\Delta t = 0.005d$) are used as reference solutions. Select different radii of the node influence domain (including $r_e = 1.001\sqrt{\Delta x^2 + \Delta y^2}$, $2.001\sqrt{\Delta x^2 + \Delta y^2}$, and $3.001\sqrt{\Delta x^2 + \Delta y^2}$) in the upwind GFDM. Fig. 4 compares the calculation results of FDM ($\Delta x = 4m$, $\Delta y = 4m$) and upwind GFDM with different radii of node influence domain on one-dimensional (1D) section $y = 40m$. It can be seen that upwind GFDM can basically achieve high accuracy, especially in pressure calculation. Only at the discontinuous water drive front, the calculation accuracy of water saturation is relatively low, but it is also satisfactory by comparing the FDM results with the same node spacing/mesh size. Fig. 5 and Fig. 6 visually compare the calculated water saturation and pressure profiles of the reference solution, FDM solution, and upwind GFDM with different radii of the node influence domain, respectively. It can be seen that the computational relative errors of various methods for pressure are small. For water saturation, the calculation results of upwind GFDM with the radius of influence domain $r_e = 1.001\sqrt{\Delta x^2 + \Delta y^2}$ is very similar with those of the upwind FDM because the GFDM discrete schemes are nearly the same as those in FDM. With the increase of the radius of the influence domain, the calculation error of upwind GFDM on water saturation decreases. This can also be seen from the water saturation profiles in Fig. 4 (b). With the increase of the radius of the node influence domain, the dissipation error at the discontinuity of the water saturation front increases. In Section 4, the detailed error analysis will be done to analyze the possible sources of the calculation error, the effects of the allocation uniformity in the node influence domain on the calculation accuracy, and the reason why the calculation accuracy of the pressure profiles does not decrease with the increase of the radius of the influence domain. The comparison results show that the upwind GFDM can almost achieve the calculation effect of the upwind FDM under the same node spacing/mesh size. However, the presented upwind GFDM is a meshless method, which can be easily applied to the calculation domain with complex geometry, while the actual underground reservoir where porous flow happens often has complex geometric characteristics (faults, complex boundaries, fractures, etc.). So, the comparisons show the great application potential of the upwind GFDM in porous flow problems.

Considering that a general numerical simulator often requires a fully implicit scheme of governing equations, this paper adopts a fully implicit nonlinear solver based on Newton iteration and automatic differentiation. Table 2 compares the Newton iteration steps under different radii of the node influence domain. It can be seen that the number of Newton iterations under different radii of the node influence domain are nearly the same, that is, the increase of the radius of the influence domain does not increase the nonlinearity of the system, indicating that the radius of the node influence domains within limits can be flexibly selected in the actual calculation. In Section 4, detailed error analysis will be done, including figuring out the sources of the calculation errors, studying the effect of the radius of the node influence domain on the computational accuracy, and roughly approximating the convergence order of the upwind GFDM. Especially, we will explain that for the studied two-phase porous flow problem, a small radius of the node influence domain should be selected to ensure enough calculation accuracy. In all, the results of this example provide a basis for constructing a robust and general numerical simulator of porous flow based on the presented upwind GFDM.

Table 1 Physical parameters used in example 1

| Parameter | Value | Parameter | Value |
| --- | --- | --- | --- |
| Permeability | 100 mD | Initial water saturation | 0.2 |
| Initial porosity | 0.3 | Irreducible water saturation | 0.2 |
| Rock compressibility | 0 MPa$^{-1}$ | Residual Oil Saturation | 0.2 |
| Oil viscosity | 10 mPa·s | Tolerated max time step size | 2 days |
| Water viscosity | 2 mPa·s | Initial time step size | 0.01 day |
| Initial pressure | 10 MPa | Tolerated error | 10$^{-6}$ |

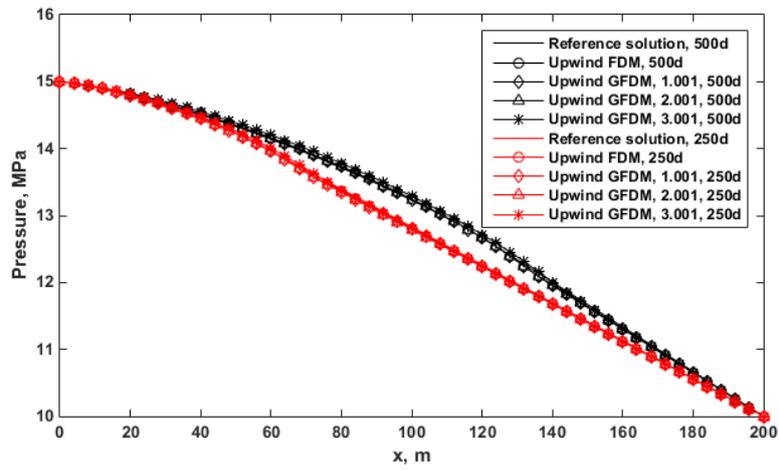

(a) pressure

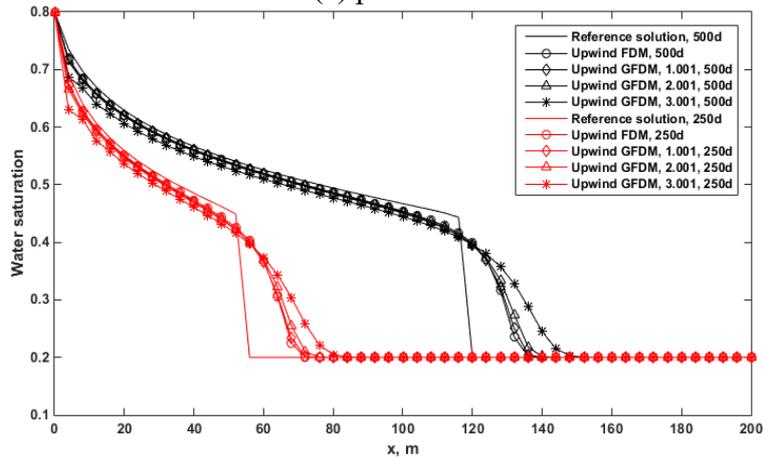

(b) water saturation

Fig. 4 Comparison of calculation results of oil-phase pressure and water saturation on $y = 40$m section

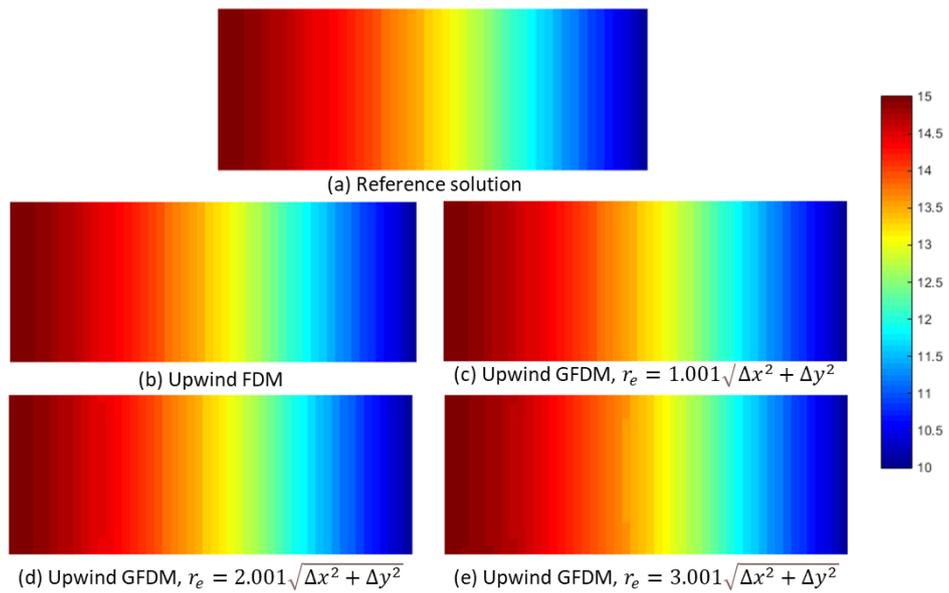

Fig. 5 The comparison of the calculated pressure profiles

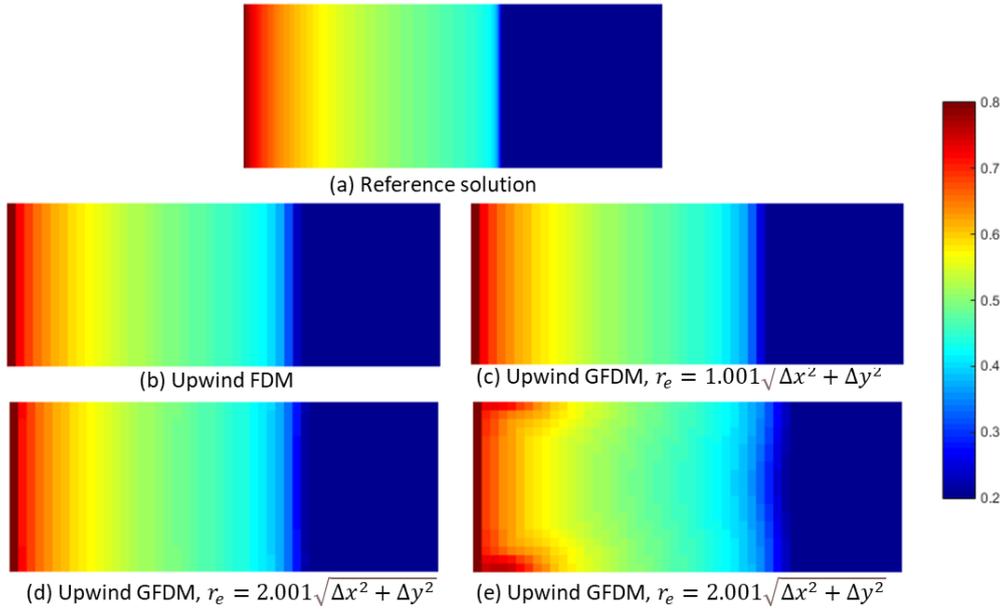

Fig. 6 The comparison of the calculated water saturation profiles

Table 2 Newton iterations in case of different radii of node influential domain

| Radius of influential domain | $1.001\sqrt{\Delta x^2 + \Delta y^2}$ | $2.001\sqrt{\Delta x^2 + \Delta y^2}$ | $3.001\sqrt{\Delta x^2 + \Delta y^2}$ |
|---|---|---|---|
| Newton iterations | 763 | 763 | 762 |

3.2 A case with irregular node clouds

Section 3.1 demonstrates that upwind GFDM can obtain calculation results with satisfactory accuracy when a regular Cartesian node cloud is used. As stated above, the significant advantage of this meshless GFDM over mesh-based methods lies in the easy discretization of the calculation domain with complex geometry. For the domain with complex geometry, irregular node clouds are generally required, so this section will test the computational performances of the GFDM based on irregular node clouds.

As shown in Fig. 7 (a), for the computational domain in Section 3.1, this section first adopts a more flexible node collocation instead of the regular Cartesian node collocation. Fig. 8 shows the calculated pressure and water saturation profiles at the 500$^{th}$ day with the radius of the node influence domain set as 8m. It can be seen that the calculation results basically match the results in Fig. 6. Then, this section gives a computational domain with irregular boundary shape and the corresponding node cloud as shown in Fig. 7 (b). The initial and boundary value conditions and governing equations are the same as those in Section 3.1, except that the shapes of the boundary $\Gamma_1$ and $\Gamma_2$ are different from those in Section 3.1. Fig. 9 shows the calculated water saturation and pressure profiles at the 300$^{th}$ day by the presented GFDM, to verify the calculation results, Cartesian mesh based FDM solutions with making the nodes outside the irregular domain to meet Dirichlet boundary conditions are used. The comparison shows that the presented upwind GFDM can effectively calculate the two-phase porous flow through irregular node collocations to discretize the computational domain with complex geometry.

It should be noted that due to the irregular node collocation in the calculation domain in this section, the method of plotting profiles adopted here is: firstly, based on the lattice points with an equal interval of 1m, interpolate the calculated node values to obtain the values on each grid point, and the values at the lattice points not in the calculation domain are assigned NAN, then the profiles can be plotted according to the values on the lattice points of Cartesian distribution.

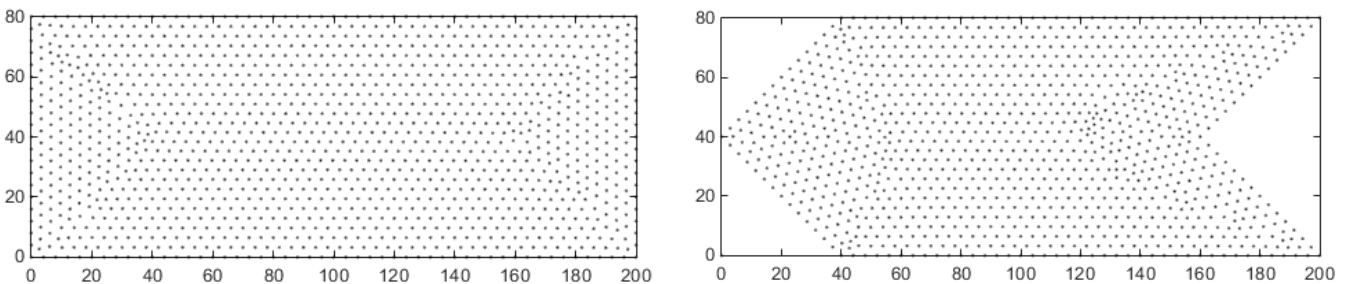

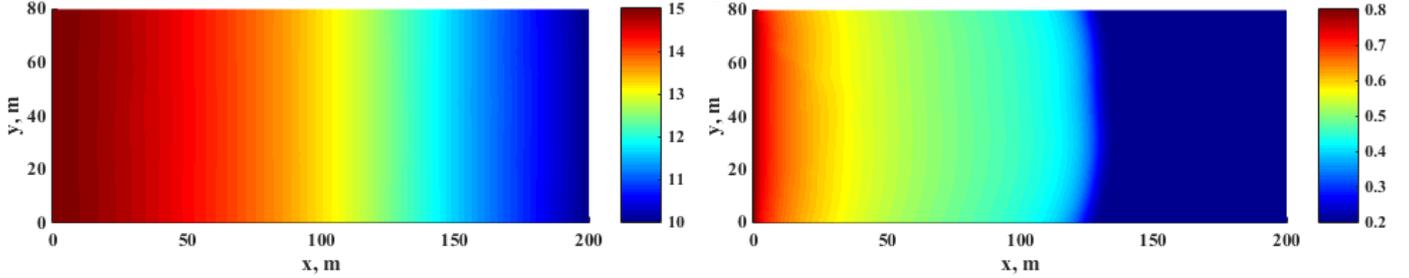

(a) rectangular domain      (b) irregular domain

Fig. 7 sketch of computational domain

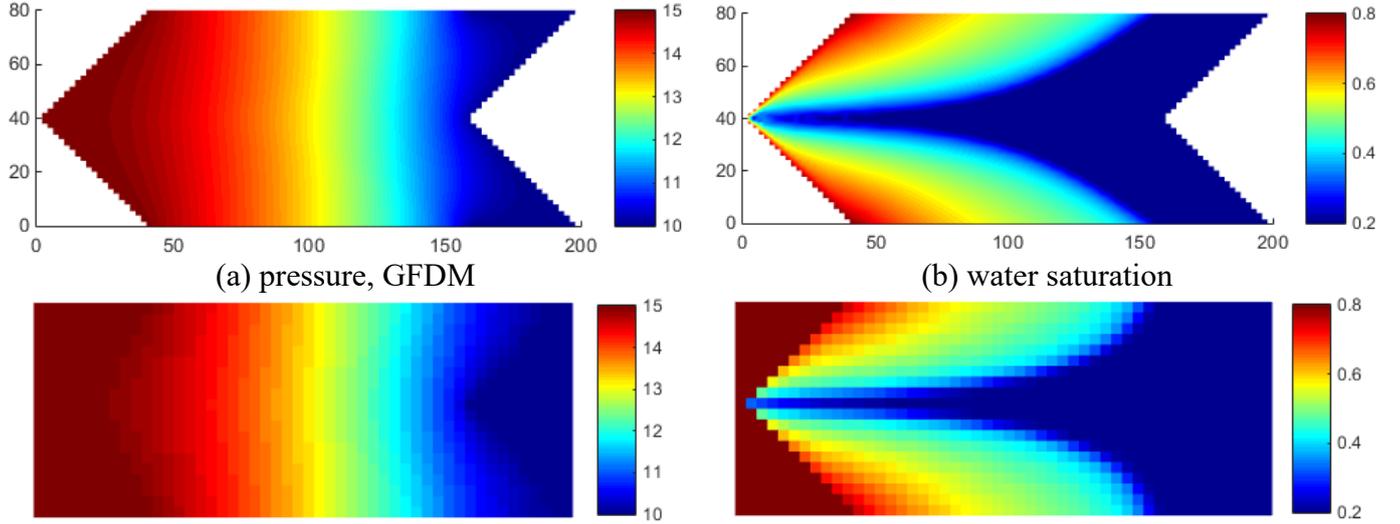

(a) pressure      (b) water saturation

Fig. 8 calculation results in case of irregular node collocation

(a) pressure, GFDM      (b) water saturation

(c) pressure, FDM      (d) water saturation, FDM

Fig. 9 Comparisons of the calculation results at the 250$^{th}$ day in case of irregular boundary shape

## 4 Error analysis
### 4.1 Effects of the allocation quality in the node influence domain on the calculation accuracy

Although the numerical example in Section 3.1 is in the two-dimensional (2D) rectangular calculation domain, it is essentially a one-dimensional (1D) water flooding case after combining the boundary conditions. Therefore, the partial derivative of the unknown functions (including pressure and water saturation) to y should be equal to zero, that is, $u_y = 0$, $u = S_w, p$.

Suppose $\Delta x = 1$, $\Delta y = 1$, when the radius of the node influence domain is set to 1.5, Fig. 10 (a) shows the node distribution in the influence domain (i.e. the gray area in the figure) of a node on the boundary (i.e. node 3) without adding virtual nodes, Fig. 10 (b) shows the node distribution in the influence domain of node 3 after adding a virtual node (i.e. the blue nodes in the figure) at each boundary node.

According to Eq, (9), the approximation expressions of the first-order derivative of the unknown function at node 3 to y (denoted as $u_{y,3}$) in Eq. (10) corresponding to Fig. 10 (a) and Fig. 10 (b) are calculated in Eq. (21) and Eq. (22) respectively. It can be seen that the difference coefficients corresponding to node 2, node 3, node 7 and node 9 in Eq. (21) are all of the order of the negative fifth power of ten, so it has little impact on the calculation results. Therefore, Eq. (22) essentially uses the first-order difference between node 3 and node 8 to approximate $u_{y,3}$. Theoretically, the difference coefficient should be equal to 1, but the difference coefficient in Eq. (21) is only 0.20988, which shows that the approximation accuracy of the partial derivative is not high without adding virtual nodes. As can be seen from Eq. (22), after adding virtual nodes, it is essentially equivalent to approximating $u_{y,3}$ with the central difference scheme between virtual node V2 and node 8. The calculated difference coefficient is very close to the theoretical value of 0.5 and has the second-order accuracy of the central difference scheme. The comparison shows that adding virtual nodes can effectively improve the approximation accuracy of the GFDM at boundary nodes. It is because of the small radius of the influence domain that the GFDM after adding virtual nodes realizes the calculation results with good accuracy like those in Fig. 6 (c). However, we should still notice that when the small radius of the node influence domain is used the calculated discontinuous water drive front like that in Fig. 6 (c) is not as clear as the finite difference solution in Fig. 6 (b), and this error will be analyzed in Section 4.2.

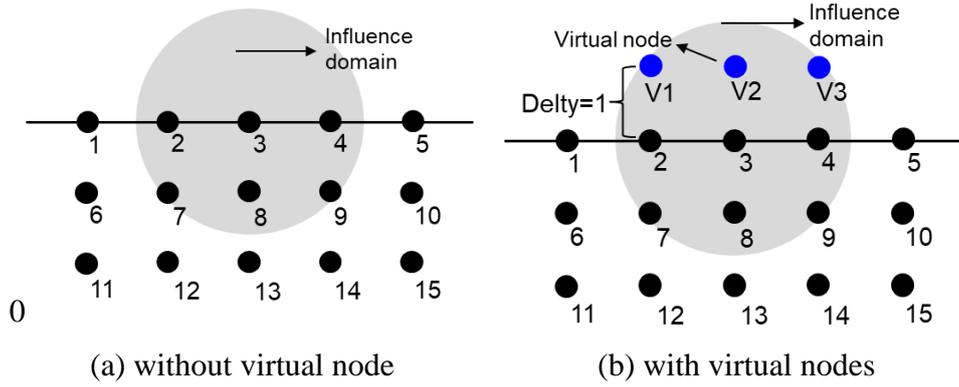

(a) without virtual node    (b) with virtual nodes

Fig. 10 sketches of the allocation nodes in the influence domain of the boundary node 3

$$u_{y,3} = 4.1537 \times 10^{-5}(u_2 - u_3) + 4.1537 \times 10^{-5}(u_4 - u_3) + 2.8212 \times 10^{-5}(u_7 - u_3) + 2.8212 \times 10^{-5}(u_9 - u_3)$$
$$+ 2.0988 \times 10^{-1}(u_8 - u_3) \tag{20}$$

$$u_{y,3} = 2.0769 \times 10^{-5}(u_2 - u_3) + 2.0769 \times 10^{-5}(u_4 - u_3) - 2.0769 \times 10^{-5}(u_7 - u_3) - 2.0769 \times 10^{-5}(u_9 - u_3)$$
$$- 4.9996 \times 10^{-1}(u_8 - u_3) + 4.9996 \times 10^{-1}(u_{V2} - u_3) \tag{21}$$

When the radius of the node influence domain increases to 2.5, Fig. 11 (a) and Fig. 11 (b) are two different ways to add virtual nodes. Fig. 11 (c) shows the node distribution in the influence domain of boundary node 3 with no virtual node added. Denote the line passing through node 3 parallel to the x-axis as $l_{3,x}$ and the line passing through node 3 parallel to the y-axis as $l_{3,y}$. It can be seen that the node distribution in the influence domain of node 3 in Fig. 11 (a) is completely symmetrical about $l_{3,x}$ and $l_{3,y}$, and Fig. 11 (b) and Fig. 11 (c) are only complete symmetry about $l_{3,y}$. Fig. 11 (c) is complete asymmetry about $l_{3,x}$, because there are no nodes above $l_{3,x}$ in Fig. 11 (c), but Fig. 11 (b) can be considered to be partially symmetrical about $l_{3,x}$, because a row of virtual nodes added above $l_{3,x}$ and a row of nodes on the lower side of $l_{3,x}$ are symmetrical about $l_{3,x}$. Next, we will analyze the influence of the symmetry or uniformity of the node distribution in the node influence domain on the calculation accuracy in detail.

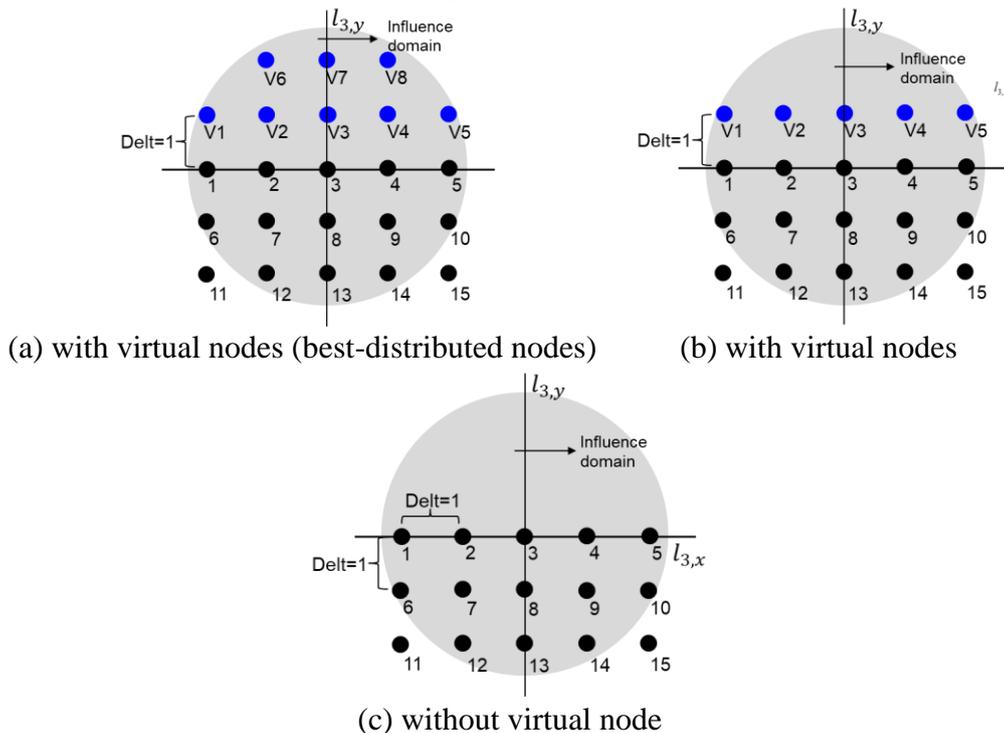

(a) with virtual nodes (best-distributed nodes)    (b) with virtual nodes

(c) without virtual node

Fig. 11 Three node distributions in the influence domain of node 3 when the radius of the influence domain

Firstly, according to the node distribution in the node influence domain shown in Fig. 11 (a), the approximation formula of $u_{y,3}$ in Eq. (23) is calculated by using Eq. (9) and Eq. (10) in Section 2.2. In the formula, the difference coefficients of nodes 1, 2, 4, and 5 are all 0 because their absolute values are less than $10^{-19}$, so they are regarded as 0.

$$\begin{aligned}u_{y,3} =\ & 0(u_1-u_3)+0(u_5-u_3)+0(u_2-u_3)+0(u_4-u_3) \\ & -2.8756\times10^{-5}(u_6-u_3)-2.8756\times10^{-5}(u_{10}-u_3)+2.8756\times10^{-5}(u_{V1}-u_3)+2.8756\times10^{-5}(u_{V5}-u_3) \\ & -7.4740\times10^{-2}(u_7-u_3)-7.4740\times10^{-2}(u_9-u_3)+7.4740\times10^{-2}(u_{V2}-u_3)+7.4740\times10^{-2}(u_{V4}-u_3) \\ & -5.7512\times10^{-5}(u_{12}-u_3)-5.7512\times10^{-5}(u_{14}-u_3)+5.7512\times10^{-5}(u_{V6}-u_3)+5.7512\times10^{-5}(u_{V8}-u_3) \\ & -0.3457(u_8-u_3)+0.3457(u_{V3}-u_3) \\ & -2.2652\times10^{-3}(u_{13}-u_3)+2.2652\times10^{-3}(u_{V7}-u_3)\end{aligned} \quad (22)$$

Since node 7 and node V2 are symmetrical about $l_{3,x}$, node 7 and node 9 are symmetrical about $l_{3,y}$, and the node distribution in the influence domain of node 3 is completely symmetrical about both of $l_{3,x}$ and $l_{3,y}$, in Eq. (23), the difference coefficient of node 7 $D_7$ is exactly the same as the difference coefficient of node 9 $D_9$. $D_7$ and the difference coefficient of node V2 $D_{V2}$ are opposite to each other. The difference coefficients $D_{V4}$ of node V4 symmetrical to the three nodes are identical with $D_{V2}$, but opposite to $D_7$ or $D_9$. To express clearly and reflect the symmetry characteristics, we denote these difference coefficients as $D_{7,9}=D_7=D_9$ and $D_{V2,V4}=D_{V2}=D_{V4}$ respectively. The difference coefficients of other nodes that are symmetrical about $l_{3,y}$ can be expressed similarly, such as $D_{2,4}=D_2=D_4$, $D_{12,14}=D_{12}=D_{14}$, $D_{V6,V8}=D_{V6}=D_{V8}$, etc.

In the numerical example in Section 3.1, the water drive direction is from left to right. When the water drive front reaches node 3, the water saturation of the node on the left of $l_{3,y}$ (including nodes 1, 2, 6, 7, 12, V1, V2, and V6) will be higher than that of the nodes on $l_{3,y}$ (including nodes 3, 8, 13, V3 and V7) and the node on the right of $l_{3,y}$ (including nodes 4, 5, 9, 10, 14, V4, V5 and V8). Since the water drive front has just reached node 3, the water saturation of the nodes on $l_{3,y}$ and the nodes on the right side of $l_{3,y}$ are almost the initial water saturation value, so the water saturation difference between the nodes on $l_{3,y}$ and the nodes on the right side of $l_{3,y}$ will be very small, while the water saturation difference between the nodes on the left side of $l_{3,y}$ and the nodes on $l_{3,y}$ will be very large, Taking nodes 2, 7, 12, V2 and V6 with the same water saturation value in theory and nodes 4, 9, 14, V4 and V8 respectively symmetrical to these five nodes about $l_{3,y}$ as an example, it can be obtained that:

$$\begin{aligned}& u_2-u_3 = u_7-u_3 = u_{12}-u_3 = u_{V2}-u_3 = u_{V6}-u_3 \\ & \gg u_4-u_3 = u_9-u_3 = u_{14}-u_3 = u_{V4}-u_3 = u_{V8}-u_3 \approx 0\end{aligned} \quad (23)$$

Then,

$$\begin{aligned}& D_2(u_2-u_3)+D_4(u_4-u_3)+D_7(u_7-u_3)+D_9(u_9-u_3)+D_{12}(u_{12}-u_3)+D_{14}(u_{14}-u_3) \\ & +D_{V2}(u_{V2}-u_3)+D_{V4}(u_{V4}-u_3)+D_{V6}(u_{V6}-u_3)+D_{V8}(u_{V8}-u_3) \\ & = D_{2,4}(u_2-u_3)+D_{2,4}(u_4-u_3)+D_{7,9}(u_7-u_3)+D_{7,9}(u_9-u_3)+D_{12,14}(u_{12}-u_3)+D_{12,14}(u_{14}-u_3) \\ & +D_{V2,V4}(u_{V2}-u_3)+D_{V2,V4}(u_{V4}-u_3)+D_{V6,V8}(u_{V6}-u_3)+D_{V6,V8}(u_{V8}-u_3) \\ & \approx (D_{2,4}+D_{7,9}+D_{12,14}+D_{V2,V4}+D_{V6,V8})(u_2-u_3)\end{aligned} \quad (24)$$

Eq. (24) reflects the strong asymmetry of water saturation distribution on both sides of $l_{3,y}$. Since $u_{y,3}$ should be equal to 0 in theory, Eq. (25) shows that it is necessary for $D_{2,4}+D_{7,9}+D_{12,14}+D_{V2,V4}+D_{V6,V8}$ to tend to 0 as much as possible to improve the approximation accuracy of Eq. (23). For the node distribution in the influence domain of node 3 which is completely symmetrical about $l_{3,x}$ and $l_{3,y}$ in Fig. 11 (a), Eq. (26) holds:

$$\begin{aligned}D_{2,4}+D_{7,9}+D_{12,14}+D_{V2,V4}+D_{V6,V8} & = D_{2,4}+(D_{7,9}+D_{V2,V4})+(D_{12,14}+D_{V6,V8}) \\ & = 0-7.4740\times10^{-2}+7.4740\times10^{-2}-5.7512\times10^{-5}+5.7512\times10^{-5}=0\end{aligned} \quad (25)$$

So, it is obtained that,

$$\begin{aligned}& D_2(u_2-u_3)+D_4(u_4-u_3)+D_7(u_7-u_3)+D_9(u_9-u_3)+D_{12}(u_{12}-u_3)+D_{14}(u_{14}-u_3) \\ & +D_{V2}(u_{V2}-u_3)+D_{V4}(u_{V4}-u_3)+D_{V6}(u_{V6}-u_3)+D_{V8}(u_{V8}-u_3)=0\end{aligned} \quad (26)$$

Similarly, it can be obtained that each node with the same water saturation value as node 6 and their symmetrical nodes meet:

$$D_1(u_1-u_3)+D_5(u_5-u_3)+D_6(u_6-u_3)+D_{10}(u_{10}-u_3)+D_{V1}(u_{V1}-u_3)+D_{V5}(u_{V5}-u_3)=0 \qquad (27)$$

Therefore, when the node distribution in the influence domain of node 3 is as shown in Fig. 11 (a), if $u$ denotes the water saturation, when the exact solution is brought into Eq. (23) that approximates $u_{y,3}$, Eq. (23) is strictly true, so the error of Eq. (23) for approximating $u_{y,3}$ is 0.

When the method of adding virtual nodes shown in Fig. 11 (b) is adopted, Eq. (29) is the calculated difference approximate expression of $u_{y,3}$.

$$\begin{aligned}u_{y,3} =& -1.3136\times10^{-5}(u_1-u_3)-1.3136\times10^{-5}(u_5-u_3)-1.0024\times10^{-3}(u_2-u_3)-1.0024\times10^{-3}(u_4-u_3)\\ &-2.8933\times10^{-5}(u_6-u_3)-2.8933\times10^{-5}(u_{10}-u_3)+2.8860\times10^{-5}(u_{V1}-u_3)+2.8860\times10^{-5}(u_{V5}-u_3)\\ &-7.4549\times10^{-2}(u_7-u_3)-7.4549\times10^{-2}(u_9-u_3)+7.5661\times10^{-2}(u_{V2}-u_3)+7.5661\times10^{-2}(u_{V4}-u_3)\\ &-5.6687\times10^{-5}(u_{12}-u_3)-5.6687\times10^{-5}(u_{14}-u_3)\\ &+0.3510(u_{V3}-u_3)-0.3438(u_8-u_3)-2.2295\times10^{-3}(u_{13}-u_3)\end{aligned} \qquad (28)$$

It can be seen that since the node distribution in the influence domain of node 3 is symmetrical about $l_{3,y}$, the difference coefficients of node 7 and node 9 symmetrical about $l_{3,y}$ are still the same, both $-7.4549\times10^{-2}$. The difference coefficients of node V2 and node V4 are also the same, both $7.5661\times10^{-2}$. Therefore, it can still be denoted as $D_{7,9}$, $D_{V2,V4}$ and $D_{2,4}$ respectively. However, since the node distribution in the influence domain is not completely symmetrical about $l_{3,x}$, although node 7 and node V2 are symmetrical about the $l_{3,x}$, their difference coefficients are not the same. It is also worth emphasizing that since the node distribution is not completely symmetrical about $l_{3,x}$ in Fig. 11 (b), the difference coefficient of node 2 on $l_{3,x}$ with the same water saturation value as node 7 in theory is no longer 0 in Eq. (23), but $-1.0024\times10^{-3}$ in Eq. (29). In this case, Eq. (30) related to nodes 2, 7, 12, V2 and nodes 4, 9, 14, V4 respectively symmetrical to these four nodes about $l_{3,y}$ holds.

$$\begin{aligned}D_{2,4}+D_{7,9}+D_{12,14}+D_{V2,V4} =& (D_{2,4}+D_{12,14})+(D_{7,9}+D_{V2,V4})\\ &-1.0024\times10^{-3}-5.6687\times10^{-5}-7.4549\times10^{-2}+7.5661\times10^{-2}=5.29\times10^{-5}\neq0\end{aligned} \qquad (29)$$

Then,

$$\begin{aligned}&D_2(u_2-u_3)+D_4(u_4-u_3)+D_7(u_7-u_3)+D_9(u_9-u_3)+D_{12}(u_{12}-u_3)+D_{14}(u_{14}-u_3)\\ &+D_{V2}(u_{V2}-u_3)+D_{V4}(u_{V4}-u_3)\\ &=D_{2,4}(u_2-u_3)+D_{2,4}(u_4-u_3)+D_{7,9}(u_7-u_3)+D_{7,9}(u_9-u_3)\\ &+D_{12,14}(u_{12}-u_3)+D_{12,14}(u_{14}-u_3)+D_{V2,V4}(u_{V2}-u_3)+D_{V2,V4}(u_{V4}-u_3)\\ &\approx(D_{2,4}+D_{7,9}+D_{V2,V4}+D_{12,14})(u_7-u_3)=5.29\times10^{-5}(u_7-u_3)\neq0\end{aligned} \qquad (30)$$

Similarly, Eq. (32) about each node with the same water saturation value as node 6 and their symmetrical nodes can be obtained.

$$\begin{aligned}&D_{1,5}(u_1-u_3)+D_{1,5}(u_5-u_3)+D_{6,10}(u_6-u_3)+D_{6,10}(u_{10}-u_3)+D_{V1,V5}(u_{V1}-u_3)+D_{V1,V5}(u_{V5}-u_3)\\ &\approx(D_{1,5}+D_{6,10}+D_{V1,V5})(u_1-u_3)=-1.3136\times10^{-5}-2.8933\times10^{-5}+2.8860\times10^{-5}=-1.3209\times10^{-5}\neq0\end{aligned} \qquad (31)$$

Therefore, it can be known that there is an error when the exact solution is brought into Eq. (29), and the absolute value of the error is positively correlated with $5.29\times10^{-5}$ and $-1.3209\times10^{-5}$.

As can be seen from Eq. (26) and Eq. (30), when approximating $u_{y,3}$, the added virtual nodes can be regarded to take the function of offsetting the error caused by the sum of the difference coefficients of the nodes with theoretically the same water saturation value.

When no virtual node is added, that is, the node distribution in the influence domain of node 3 shown in Fig. 5 (c), Eq. (33) is the calculated difference approximate expression of $u_{y,3}$.

$$\begin{aligned}u_{y,3} =& 3.1575\times10^{-3}(u_1-u_3)+3.1575\times10^{-3}(u_5-u_3)+2.4093\times10^{-1}(u_2-u_3)+2.4093\times10^{-1}(u_4-u_3)\\ &-4.2024\times10^{-5}(u_6-u_3)-4.2024\times10^{-5}(u_{10}-u_3)-2.6549\times10^{-1}(u_7-u_3)-2.6549\times10^{-1}(u_9-u_3)\\ &+1.2100\times10^{-2}(u_{12}-u_3)+1.2100\times10^{-2}(u_{14}-u_3)\\ &-1.4689(u_8-u_3)+0.4758(u_{13}-u_3)\end{aligned} \qquad (32)$$

It can be seen that in this case, since the node distribution in the influence domain of node 3 is still symmetrical about $l_{3,y}$, the difference coefficients of node 7 and node 9 are still the same. but the node distribution is completely asymmetric about $l_{3,x}$. At this time, Eq. (34) subject to nodes 2, 7, 12 with the same water saturation value in theory and corresponding nodes 4, 9, and 14 respectively symmetrical about $l_{3,y}$ can be obtained.

$$D_{2,4} + D_{7,9} + D_{12,14} = 2.4093 \times 10^{-1} - 2.6549 \times 10^{-1} + 1.2100 \times 10^{-2} = -1.25 \times 10^{-2} \neq 0 \tag{33}$$

Then,

$$\begin{aligned} &D_2(u_2 - u_3) + D_4(u_4 - u_3) + D_7(u_7 - u_3) + D_9(u_9 - u_3) + D_{12}(u_{12} - u_3) + D_{14}(u_{14} - u_3) \\ &= D_{2,4}(u_2 - u_3) + D_{2,4}(u_4 - u_3) + D_{7,9}(u_7 - u_3) + D_{7,9}(u_9 - u_3) + D_{12,14}(u_{12} - u_3) + D_{12,14}(u_{14} - u_3) \\ &\approx (D_{7,9} + D_{2,4} + D_{12,14})(u_7 - u_3) = -1.25 \times 10^{-2}(u_7 - u_3) \neq 0 \end{aligned} \tag{34}$$

It can be seen that the absolute value of the error coefficient $-1.25 \times 10^{-2}$ in Eq. (34) is greater than the absolute value of the error coefficient $5.29 \times 10^{-5}$ in Eq. (30), and much significantly greater than the error coefficient 0 in Eq. (26). It shows that with the decrease of the symmetry of the node distribution in the influence domain of node 3, the ability to offset the error between nodes in the node influence domain gradually decreases, resulting in the gradual decrease of the approximation accuracy of the derivatives of the unknown function in GFDM.

Of course, strict symmetry generally occurs in the case of Cartesian collocation, while for general collocation, symmetry corresponds to the uniformity of node distribution in the node influence domain, that is, there are almost as many and symmetrical nodes as possible on both sides of a straight line passing through the central node, this is an important factor affecting the calculation accuracy of GFDM based on collocation. The closer the gravity center of the node distribution in the influence domain of a node is to the node, the better the uniformity of the distribution point may be considered

The comparison between the results in Fig. 4 and Fig. 6 in Section 3.1 also shows that the increase of the radius of the node influence domain leads to the decline of calculation accuracy. This is because to build a stable GFDM based simulator, in this work a virtual node is added correspondingly to each boundary node along the external normal vector at the boundary node with a certain distance. For the Cartesian collocation in the example, the distance is just $\Delta x$ or $\Delta y$. That is to say, the sketch of adding the virtual nodes is the same as Fig. 11 (b). Therefore, with the gradual increase of the influence domain, the symmetry or uniformity of the node distribution in the influence domain of the boundary nodes and the nodes close to the boundary becomes worse and worse, resulting in the gradual decline of the accuracy of the GFDM based calculation for the two-phase porous flow problem in this paper. Moreover, it can be predicted that because the node distribution symmetry or uniformity of the nodes near the boundary and the internal nodes are different, for the nodes with the same *x* coordinate, the error of estimating the partial derivative of water saturation to *y* at these nodes by using the GFDM theory will be different, resulting in the bending of the water drive front in theory. The bending water drive fronts shown in Fig. 6 (d) and Fig. 6 (e), and the more bending water drive front in Fig. 6 (e) corresponding to the larger radius of the node influence domain, verify the analysis in this section.

4.2 Flow physics-based analysis

The analysis in Section 4.1 demonstrates the influence of the collocation quality (i.e. symmetry or uniformity) in the node influence domain on the calculation accuracy of the hyperbolic two-phase porous flow problem studied in this paper. This section will analyze another error from the perspective of the actual physics of two-phase porous flow. When the radius of the influence domain is larger, the calculation results of GFDM will become more unphysical, and the discontinuity reflected in the water drive front will become less and less obvious, which will reduce the calculation accuracy.

As shown in Fig. 12, the direction of water flooding in the example in Section 3.1 is from left to right. Suppose the node distribution of the calculation domain is shown in Fig. 12, and node 1 is a node on the left boundary and node 2 is a node on the right boundary. At the initial time *t*, the water saturation at node 2 is irreducible water saturation. Imagine an extreme case. In this case, the radius of the influence domain of node 1 is large enough to include node 2. No matter what the time step of the Newton iteration convergence in nonlinear solution is, it might as well be set as a small time step (because the smaller the time step, the easier it is to converge). Since node 2 is within the influence domain of node 1, the water-phase transfer quantity between node 1 and node 2 can be measured by

$$\frac{k_{12}^{t+\Delta t} k_{rw}\left(S_{w,1}^{t+\Delta t}\right)}{\mu_{w,12}^{t+\Delta t}} \left(e_{3,j}^i + e_{4,j}^i\right)\left(p_{w,1}^{t+\Delta t} - p_{w,2}^{t+\Delta t}\right)\Delta t \tag{35}$$

and $k_{rw}\left(S_{w,1}^{t+\Delta t}\right) > 0$ and $\left(p_{w,1}^{t+\Delta t} - p_{w,2}^{t+\Delta t}\right) > 0$ hold, so the water saturation at node 2 will rise at the time $t + \Delta t$, which means that the water flowing from the left boundary reaches the right boundary in a very small time step $\Delta t$, which is obviously unphysical. In other words, the larger the radius of the node influence domain, the greater the dissipation error, and the less obvious the discontinuity of the water drive front. In Fig. 4 (b) and Fig. 6 (c), (d) and (e), the larger the radius of the node influence area, the less obvious the discontinuity at the water drive front, which verifies the analysis in this section.

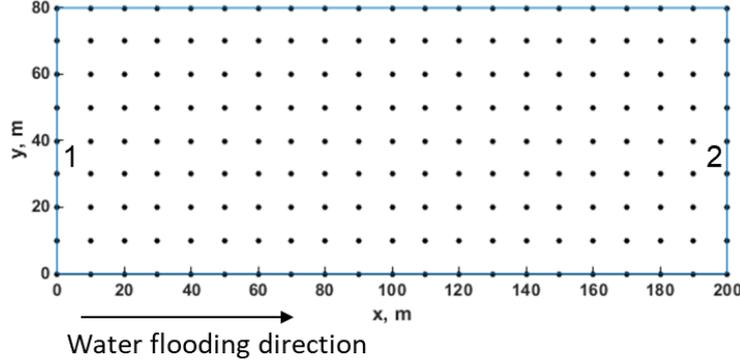

Fig. 12 the sketch of the extreme case

4.3 The Difference between GFDM based analysis of the hyperbolic two-phase porous flow problem and elliptic problems

In Sections 4.1 and 4.2, two sources of the calculation error of the water saturation are analyzed, including the effect of the node collocation quality in the node influence domain on the accuracy of the generalized difference operator analyzed in Section 4.1 and the unphysical error analyzed in Section 4.2. However, it should also be noted that the calculation errors of pressure distribution in Fig. 4 (a) and Fig. 5 in Section 3.1 under different radii of the node influence domain are very small. Since the rock compressibility coefficient in the example in Section 3.1 is zero (shown in Table 1), the pressure function actually satisfies the following Eq. (37) by adding Eq. (1) and Eq. (2), that is, it is a nonlinear elliptic equation. Considering that the size of the node influence region is often large and can obtain high enough calculation accuracy in previous studies of GFDM solutions of elliptic equations, we guess that the hyperbolic two-phase porous flow problem (especially the water saturation) concerned in this paper have different requirements for the radius of the node influence region from the elliptic problems when solved by GFDM, The elliptic problem tends not to limit the radius of the node influence domain, which should be a significant difference between the studied hyperbolic two-phase porous flow and the elliptic problems.

$$\alpha \nabla \cdot \left( \left( \frac{kk_{ro}}{\mu_o} + \frac{kk_{rw}}{\mu_w} \right) \nabla p_o \right) + q_o = \frac{\partial \phi}{\partial t} = 0 \tag{36}$$

Next, we will demonstrate the above conjecture by considering that $u$ denotes the pressure function based on the analysis in Section 4.1. Taking node 7 and node 9 that are symmetrical about $l_{3,y}$ as an example, the analysis in Section 4.1 is based on the result that when $u$ denotes the water saturation, there will be $u_7 - u_3 \gg u_3 - u_9$, which leads to the consideration of the symmetry or uniformity of the node distribution in the node influence domain. However, when only for elliptic problems, such as the pressure function satisfying Laplace equation, Poisson equation or the nonlinear elliptic equation in Eq. (37), Eq. (38) generally holds.

$$u_7 - u_3 \approx u_3 - u_9 \tag{37}$$

When Eq. (38) holds, still take Fig. 11 as an example, since the node collocations in the influence domain of nodes 3 in Fig. 11 (a), (b), and (c) are all symmetrical about $l_{3,y}$, when approximating $u_{y,3}$, the difference coefficient of node 7 $D_7$ must be equal to the difference coefficient of node 9 $D_9$, which can also be uniformly denoted as $D_{7,9}$, then there is

$$D_7(u_7 - u_3) + D_9(u_9 - u_3) = D_{7,9}\left[(u_7 - u_3) - (u_3 - u_9)\right] \approx 0 \tag{38}$$

Similarly, for any node $i$, let its node symmetrical about $l_{3,y}$ be node $j$, then Eq. (40) holds.

$$D_i(u_i - u_3) + D_j(u_j - u_3) = D_{i,j}\left[(u_i - u_3) - (u_3 - u_j)\right] \approx 0 \tag{39}$$

In this case, using either Eq. (23) corresponding to Fig. 11 (a), Eq. (29) corresponding to Fig. 11 (b) or Eq. (33) corresponding to Fig. 11 (c), the calculated $u_{y,3}$ tend to 0 with much small error. Therefore, for the unknown function satisfying the elliptic equation, the radius of the node influence region has little effect on the calculation accuracy.

The above analysis demonstrates that when the GFDM is applied to elliptic problems, there is no significant requirement of the symmetry or uniformity of the node collocation in the node influence domain, which is why the radius of the node influence domain can be large but without reducing the calculation accuracy in many previous GFDM based numerical studies of elliptic problems. The discovery of this difference on the requirement of the radius of the node influence domain between hyperbolic two-phase porous flow problem and elliptic problems may be also an important innovation of this work.

4.4 The convergence order of this upwind GFDM for the two-phase porous flow problem

According to the above analysis in Sections 4.1 and 4.2, it can be known that when the node cloud of in the calculation domain is determined, the radius of the node influence domain radius needs to be small to improve the calculation accuracy. However, if the radius of the influence domain is set too small, it will degenerate into the traditional finite difference method in the case of the Cartesian node cloud in Section 3.1. In order to explore the convergence order of the upwind GFDM when upwind GFDM is different from traditional FDM, this section set $\Delta x = \Delta y$ and equals 0.5, 1, 2, 4, 8, 10, 20, and 40 respectively, and set the radius of the node influence domain as $r_e = 1.5\Delta x$. At this time, the node in the diagonal direction of the central node has sufficient weight, which ensures the GFDM discrete schemes different from those in FDM. Fig. 13 compares the relative calculation error of the pressure and water saturation calculated by the upwind GFDM and FDM versus node spacing/mesh size, in which the relative error is defined in Eq. (40). It can be seen that in the double logarithmic coordinate system, when the mesh size is large, the slope of the relative error curve of the upwind GFDM is nearly the same as that of upwind FDM. While when the mesh size is small, the slope of the relative error curve of the upwind GFDM is obviously less than that of upwind FDM, indicating that the convergence order of the upwind GFDM is less than or equal to that of the upwind FDM. For this phenomenon, we think that the generalized difference operator based on the moving least squares generally cannot make the minimized weighted truncation error equal to 0, so that the generalized difference expressions in Eq. (10) have the error with low order. When the node spacing is large, the proportion of the low-order error is smaller than the high-order error which both GFDM and FDM have, so the convergence order of GFDM and FDM is almost the same. When the node spacing is small, the proportion of low-order error in GFDM in the total error is increasing, resulting in the convergence order of GFDM is significantly lower than that of FDM.

$$RE = \frac{\|\mathbf{u} - \mathbf{u}_{ref}\|_2}{\|\mathbf{u}_{ref}\|_2} \tag{40}$$

where $\mathbf{u}$ is the vector consists of calculated nodal values of the pressure or the water saturation. $\mathbf{u}_{ref}$ is the vector consists of high-resolution solutions calculated by the upwind FDM. $\|\cdot\|_2$ is the 2-norm of a vector.

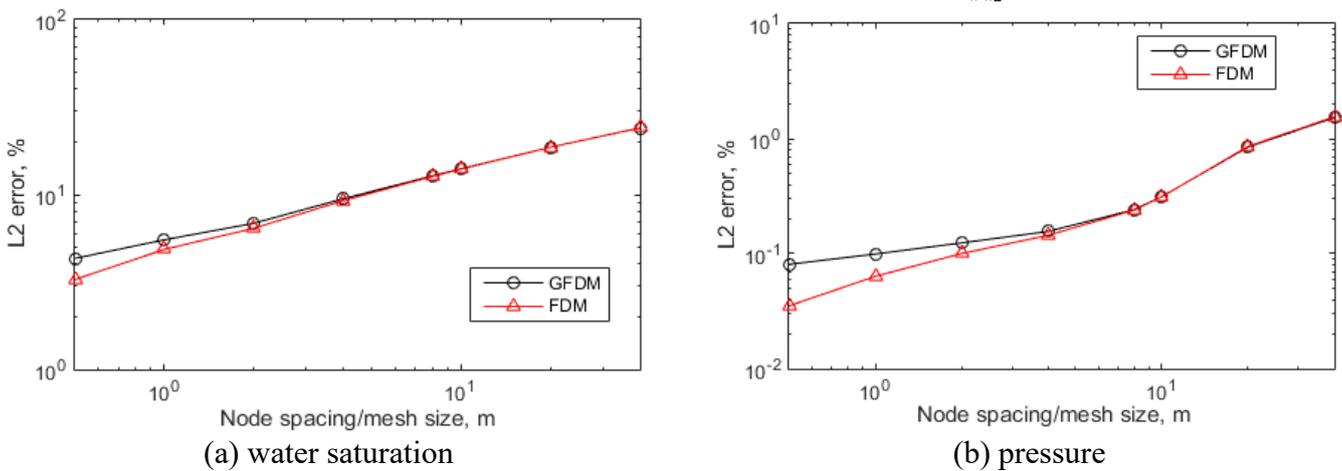

(a) water saturation  (b) pressure
Fig. 13 the relative calculation errors versus the node spacing/mesh size

5. Conclusions
Throughout the whole paper, four main conclusions can be obtained as follows:
(i) The upwind GFDM is first applied to two-phase porous flow equations, and realizes effective meshless calculation. Compared with the traditional mesh-based methods, this method uses the node cloud to discretize the computational domain, avoids the difficulty of high-quality mesh generation, and has the advantages of more flexibility and simplicity in dealing with the computational domain with complex geometry.

(ii) Numerical examples and error analysis study the effect of the allocation uniformity in the node influence domain on the calculation accuracy, and finds that with the gradual increase of the influence domain, the symmetry or uniformity of the node distribution in the influence domain of the boundary nodes and the nodes close to the boundary becomes worse, resulting in the gradual decline of the accuracy of the GFDM based calculation for the two-phase porous flow problem in this paper.

(iii) The radius of the node influence domain is required small to not only reduce the negative effect of a large radius on the allocation quality, but also reduce the unphysical dispersion error.

(iv) A significant difference between the GFDM based solutions of the studied hyperbolic two-phase porous flow problem and elliptic equations lies in the effect of the radius of the node influence domain on the calculation accuracy. For elliptic problems, the radius of the influence domain has little impact on the accuracy.

(v) GFDM based on the moving least squares generally cannot make the minimized weighted truncation error equal to 0, so that the generalized difference expressions have the error with low order. When the node spacing is large, the proportion of the low-order error is smaller than the high-order error which both GFDM and FDM have, so the convergence orders of GFDM and FDM are almost the same. When the node spacing is small, the proportion of low-order error in GFDM in the total error is increasing, resulting in the convergence order of GFDM is significantly lower than that of FDM.

6. Acknowledgements

Dr. Rao thanks the supports from the National Natural Science Foundation of China (Nos. 52104017, 51874044, 51922007), the Open Foundation of Cooperative Innovation Center of Unconventional Oil and Gas (Ministry of Education & Hubei Province) (No. UOG2022-14), and the open fund of the State Center for Research and Development of Oil Shale Exploitation (33550000-21-ZC0611-0008).

7. Author contributions
Xiang Rao: Conceptualization; Data curation; Formal analysis; Funding acquisition; Investigation; Methodology; Project administration; Software; Validation; Visualization; Writing - original draft; Writing - review & editing;
Yina Liu: Data curation; Investigation;
Hui Zhao: Investigation; Funding acquisition;